\documentclass[12pt]{amsart}
\usepackage{amssymb, amscd, amsmath, amsthm, epsf, epsfig, latexsym, enumerate}

\theoremstyle{definition}

\numberwithin{equation}{section}

\begin{document}
\title{Unknown Values in the Table of Knots}
  
\author{Jae Choon Cha}
\address{Department of Mathematics \\
Pohang University of Science and Technology\\
Pohang Gyungbuk 790--784\\
Republic of Korea}
\email{jccha@postech.ac.kr}
\author{Charles Livingston} 
\address{Department of Mathematics\\
 Indiana University\\ 
 Bloomington, IN   47405}
  \email{livingst@indiana.edu}
 
 \thanks{Research supported by  NSF grant DMS--0406934.} 
\maketitle
 \begin{abstract} This paper, to be regularly updated,  lists those prime knots with the fewest possible number of crossings for which values of basic knot invariants, such as the unknotting number or the smooth 4--genus, are unknown.  This list is being developed in conjunction with {\sf KnotInfo} ({\tt www.indiana.edu/$\tilde{\ }$knotinfo}), a web-based table of knot invariants. 
 \end{abstract}
\section{Introduction}    
 For many basic invariants of  knots in $S^3$   there is no known algorithm for their computation.  Some invariants are  known to be algorithmic, but the algorithms are not efficient for computations for knots with more than a few crossings. This paper lists the prime knots among those of few crossings  for which the value of some of these invariants is not known.  
 
 Definitions of each of these invariants are given in the sections that follow.  References in the initial version of this paper are limited, but will be included as progress is made and contributions are received.  More details can be found on the accompanying website  {\sf KnotInfo} \cite{ki} at   {\tt www.indiana.edu/$\tilde{\ }$knotinfo}, presenting an updated list of values of knot invariants for low-crossing number knots.  
 \vskip.1in
 \noindent {\bf Conventions for tabulation}  In the table that follows, if a knot is listed of $n$--crossings, then the table lists all prime knots of that number of crossings or less  for which the value is unknown.  
 Some values include entries with question marks, such as ``12?''  For these  invariants the values for knots with fewer crossings are known.  The notation ``12?'' indicates that the value is known for all prime knots of 11 crossings or less  and that the computation for 12--crossing prime knots has not been undertaken yet.    \vskip.1in
 
 \noindent{\bf  Submissions}  It is the intention of the authors to keep this paper up-to-date.  We cannot verify all the results submitted to us, so will  include new information only when we can reference a source that readers can refer to; these can be posted papers or websites with more information.  
 
 \vskip.5in 
 \begin{table}[h]
 \begin{tabular}{||  l   |   r   ||} \noalign{\hrule height  1pt}
 {\bf  Invariant } & {\bf List of knots}\\ \noalign{\hrule height  1pt}
 A-Polynomial & ($8_{5,10, 15-21}$)*, ($9_{30,32,33,34,39,40}$)* \\ \hline
 Arc Index  (Grid Number)&  $13n?$  \\ \hline
Braid Index& $13?$\\ \hline
Bridge Index& $13?$\\ \hline
 Concordance Genus (smooth)& ($19$ unknown) $11$?\\ \hline
Concordance Order (smooth)& $ 11_{n34}, 12_{n846}$\\ \hline
Concordance Order (top)&  $12a_{550,644,690,774,1139,1140}$, $12n_{846}$ \\ \hline
Crosscap Number (non-orientable 3--genus)&  $10_{157,158,159,163,164}$\\ \hline
Fibered& $13$? \\ \hline
Genus (3--genus)& 13?\\ \hline
Genus (smooth 4--genus)&$11n_{80},$ (21 12-crossing unknown)\\ \hline
Genus (topological 4--genus)& $ 12a_{244,810,905,1142}, 12n_{549,555,642}$\\ \hline
Morse--Novikov Number& $  11  $?\\ \hline
Nakanishi Index & $11?$ \\ \hline
Oszv\'ath--Szab\'o $\tau$ Invariant& $11n{?}$\\ \hline
Polygon Number& $8_{1-15}$?\\ \hline
Slice (Topological)& $ 13?$\\ \hline
Slice (Smooth) & $ 13?$\\ \hline
Super-Bridge Number& $5_2, 6_1, 6_2, 6_3$\\ \hline
Tetrahedral Number& $7_6, 7_7$\\ \hline
Thurston--Bennequin Number& $13$?  \\ \hline
Turaev Genus& $11n_{118} $($24$ $12$-crossing unknown) \\ \hline
Unknotting Number&  $10_{11, 47 ,51, 54, 61, 76, 77,79, 100} $\\ \hline
Unknotting Number 1& $11_n$?\\ \hline
 \noalign{\hrule height  1pt}
 \end{tabular}
  * See discussion below. $\hskip3in$
 \vskip.2in
 \caption{Table of Unknown Values of Knot Invariants}
 \end{table}
   \section{Definitions}
 
  \subsection{A-Polynomial}  There is a map of the  $SL_2({\bf C})$ representation space of a knot complement to ${\bf C}^* \times {\bf C}^*$, given by evaluating the eigenvalue of the  representation on the meridian and longitude. The closure of the image is a   variety  defined by a single polynomial, called the A-Polynomial.  Jim Hoste informs us that for 2--bridge knots the computation of the A-Polynomial has a fairly efficient algorithmic and has provided us with a list of values for 2--bridge knots through 9 crossings, but beyond that only sporadic values are known.
  
  {\bf Update July 7, 2006 and May 15, 2008:} There is an alternative approach to defining the A-polynomial, based on
varieties arising from the gluing equations of an ideal hyperbolic
tetrahedral decomposition of the knot complement.  These equations
represent the restrictions on the defining parameters of the ideal  
tetrahedron that arise from the requirement that the resulting space
have a hyperbolic structure. The holonomies of the meridian and longitude can also be expressed as polynomials in the tetrahedral parameters, called the completeness equations. Eliminating the
tetrahedral parameters from the gluing and completeness equations
gives a 2-variable polynomial which in general divides
the $PSL_2({\bf C})$ version of the A-polynomial, and
can depend on the choice of decomposition.  For further discussion, see
the appendix of~\cite{brd} written by Dunfield or~\cite{c}.  Marc Culler has
computed the A-polynomials associated to such tetrahedral
decompositions; the first unknown values are for $9_{30,32,33,34,39,40}$.  
 
   \subsection{Arc Index}  Every knot has an embedding in ${\bf R}^3$ so that it is the union of closed arcs, each of which lies in a distinct half plane, with the half planes meeting in single line, the common boundary of the half planes.  The minimum number of arcs required is the arc index of the knot.  We thank Peter Cromwell  for providing us with the basic data (see~\cite{cr}) and note that Bae and Park~\cite{bp} have shown that the arc index of an alternating knot is 2 more than its crossing number. The arc index of all knots with 10 or fewer crossings has been computed, and we are aware of no tabulation for 11--crossing nonalternating knots.
   
    {\bf Update January 1, 2007:} Ng~\cite{ng3} has completed the calculation of arc indices for all 11 crossing knots.
   
    {\bf Update September,  2008} Litherland has computed the values for most unknown cases of 12 crossing knots. The remaining cases are: $12n_{41, 119, 120, 121, 145, 153,  199, 200, 243,  260, 282, 310, 322, 351, 362, 368, 377, 403, 414, 425, 475, 523, 549}$.
    
      {\bf Update October,  2008} Gyo Taek Jin, et al. ~\cite{jkl}  have completed the determination of the arc indices of 12 crossing knots. 
 \subsection{Braid Index}  Every knot can be described as the closure of a braid.  The braid index of the knot is the minimum number of strands required for such a braid.
 
 \subsection{Bridge Index}  Given a regular diagram for a knot in the $(x,y)$--plane, the projection onto the line $x=0$  has a finite number of local maximum points.  The minimum number of such local maxima, taken over all possible diagrams, is the bridge index of the knot.  {\bf Update June, 2014:}  Musick~\cite{music} has resolved all remaining cases of 11-crossing knots.
 
 \subsection{Concordance Genus}  The concordance genus of a knot $K$ is the minimum 3--genus among all knots concordant to $K$.  The value of this invariant depends on the category. For instance, for a topologically slice  but not smoothly slice knot, the topological concordance genus would be 0  but the smooth concordance genus would be positive.  For the three unknown knots of crossing number 10 or less, we do not know if the distinction of category appears.  {\bf Update September 17, 2008:}  Livingston~\cite{liv1} has computed the values for these last three cases of knots with fewer than 11 crossings.
 
 \subsection{Concordance Order}  The concordance order of a knot is the order of the element it represents in the concordance group.  There is a map of the concordance group to Levine's algebraic concordance group, having elements only of order 2, 4 and infinite.  All knots of 10 crossings or less with algebraic order 4 are known to have infinite concordance order.  Some knots of algebraic order 2 are (or are concordant to) negative amphicheiral knots, and thus are of order 2.  Other than these, all those of algebraic order 2 are known to have concordance order greater than 2, with the exception of $10_{158}$, which has unknown order.   {\bf Update January 1, 2007:} Jablan and Naik have applied Heegaard-Floer homology to rule out smooth concordance order 4 for several 2--bridge knots. {\bf Update April 20, 2008:}  Grigsby, Ruberman and Strle~\cite{grs} have further applied Heegaard Floer invariants to study the order of 2--bridge knots and Lisca~\cite{lis} has completed the analysis of the smooth concordance order of 2--bridge knots. {\bf Update May 23, 2008:} Levine~\cite{le} has applied the techniques  of~\cite{grs} to resolve 46 of the remaining unknown orders for knots of 11 crossings or less, in the smooth category.

 \subsection{Crosscap Number (non-orientable 3--genus)}  The crosscap number of a knot $K$ is the minimum value of $n$ so that $K$ bounds a surface homeomorphic to  $\#_nP^2 - B^2$ in $S^3$.  This is also called the {\it non-orientable genus} of a knot.  Hirasawa and Teragaito~\cite{ht} have described an algorithm for computing the crosscap number of a 2--bridge knot.  {\bf Update July,  2011:} Burton and Ozlen~\cite{bo} have used normal surfaces to find nonorientable surfaces of low genus for prime knots of 12 or fewer crossings, providing new upper bounds for 778 of these knots.  In particular, they resolved the crosscap number of the knot $8_{20}$. {\bf Update June, 2014:} Adams and Kindred~\ref{adamskindred} have completed the determinations of crosscap numbers for 9 crossing knots.
   \subsection{Fibering}  There is no known efficient algorithm for determining if a knot is fibered.  For knots with 10 or fewer crossings the classification has been known for some time.  For 11 crossings, several authors have undertaken the determination of which are fibered; we are told that Gabai was the first to complete the determination for 11 crossing knots.  For 12 crossing knots, work of Stoimenow, done in part with  Hirasawa, determined most cases.  These results are posted at~\cite{stoi}.   Friedl and Kim~\cite{fk} have used twisted Alexander polynomials to  obstruct the fibering of the remaining cases.  Rasmussen has used knot Floer homology, defined by  Ozsv\'ath-Szab\'o, to confirm the results of Stoimenow,  Hirasawa, Kim and Friedl.
 
 \subsection{Genus (3--genus)}  The genus of a knot is the minimum genus among all embedded orientable surface in $S^3$ with boundary the knot.
  
 \subsection{Genus (smooth 4--genus)}  The smooth 4--genus  of a knot is the   minimum genus among all smoothly embedded orientable  surfaces in the 4--ball having boundary the knot.  {\bf Update January 15, 2006:} The last remaining unknown case for 10 crossing knots, that of $10_{51}$ was resolved by Selahi Durusoy, who observed that a single crossing change results in the slice knot $8_8$.  {\bf Update July 9, 2006:} Given the completion of the table for 10--crossing knots, a preliminary effort to compute the smooth and topological 4-genus of 11-- and 12--crossing knots has been undertaken.  The values for roughly half of the 11--crossing knots are unknown, and for 11 of these knots the gap in the possible value is 2.

 \subsection{Genus (topological 4--genus)}  The topological 4--genus  of a knot is the minimum genus among all topologically locally flat embedded orientable surfaces in the 4--ball having boundary the knot. {\bf Update July 9, 2006:} See the update for the smooth 4--genus, above.
 
\subsection{Morse--Novikov Number}  The projection function from the peripheral torus in a knot complement to the meridian which is constant on a longitude extends to the entire knot complement.  Among such extensions having only nondegenerate critical points, the minimum number of such critical points is the Morse Novikov Number of the knot.  The value of this invariant is 0 if and only if the knot is fibered.  Also, the invariant is always even.  Basic results are contained in a paper by Goda,~\cite{g}.
 \subsection{Nakanishi Index.}  The Alexander module of a knot has a square presentation.  The minimum size of such a presentation   is the Nakanishi index.
 
 \subsection{Oszv\'ath--Szab\'o Invariant}  Ozsv\'ath and Szab\'o have defined a homomorphism $\tau$ from the concordance group of knots to ${\bf Z}$ using their theory of Heegaard Floer homology.  Its value is called the Ozsv\'ath--Szab\'o invariant of  the knot.  For alternating knots it is determined by the classical signature.  {\bf Update January 1, 2007:} Recent work of Manolescu, Ozsv\'ath, Sarkar, Szab\'o, and D.~Thurston, has provided combinatorial algorithms for the computation of various forms of Heegaard-Floer and knot Floer homology.  These have been implemented by Baldwin and Gillam~\cite{bd} to provide the value for $10_{141}$, the last unknown 10 crossing knot.

 \subsection{Polygon Number} The polygon number of a knot is the minimum number of vertices required in a polygonal description of the knot. We thank Peter Cromwell~\cite{cr}  for providing us with the basic data (see~\cite{cr}), summarizing the results of Randell, Negami, and Calvo.  The polygon numbers of knots with 7 crossings or less is complete.  Of 8 crossing knots only $8_{16-21}$, are known.
 
 \subsection{Slice (Topological and Smooth)}  A knot is slice if its 4--ball genus is 0.  There are topologically slice knots that are not smoothly slice, but there are no known examples of this with 11 or fewer crossings, so the two categories are not distinguished in the table.  {\bf Update July 9, 2006:} Slicing disks for several 11-crossing knots ($11_{a28,a35,a96}$) have been found by McAtee. {\bf Update April 20, 2008:}  Alexander Stoimenow conducted a computer search for smooth  slice disks. This resolved the (topological) slicing question for all but 18 of the prime knots of 12 or fewer crossings.  (In the smooth setting, the knot $11n_{34}$, which has Alexander polynomial 1, has unknown smooth 4--ball genus.)  Herald, Kirk, and Livingston~\cite{hkl} used Casson-Gordon signature invariants to rule out sliceness for 2 of these 18, and developed methods of twisted polynomials and Casson-Gordon invariants to rule out another 15.  The only remaining case is the knot $12a_{631}$.  {\bf Update June 2014:} Axel Seegler has shown that $12a_{631}$ is slice.

 \subsection{Super-Bridge Number} For a given diagram of a knot in the $(x,y)$--plane,   there is a maximum value for the number of local maxima for orthogonal projections to lines in the plane, taken over all lines in the plane for which the projection map is regular.  The smallest value of the maximum, taken over all possible diagrams of a knot is called the super-bridge number.  Randell~\cite{r} has observed that the super-bridge number is bounded above by half the polygon number.
 
 \subsection{Tetrahedral Number}  The complement of a hyperbolic knot can be decomposed into ideal hyperbolic tetrahedron.  The minimal number required is the tetrahedral number of a knot.  A census of hyperbolic knots in terms of tetrahedral decompositions is now complete for those with tetrahedral number 7 or less~\cite{brd, cdw}.  Values can be found in~\cite{ki} using the naming convention choice of ``Tetrahedral Census Name."   
 
 \subsection{Thurston--Bennequin Number}  Every knot has a Legendrian representative  and the Thurston--Bennequin number of such a representative is defined using that Legendrian structure.  The Thurston--Bennequin number of a knot is the maximum value of that invariant, taken over all possible Legendrian representatives.

 There is also a combinatorial definition.  Every knot has a diagram such that at each crossing one strand intersects the bottom right and top left quadrants (formed by the vertical and horizontal lines at the crossing) and the other intersects the bottom left and top right quadrant.  Furthermore it can be arranged  that the strand that meets the bottom right quadrant passes over the other strand.  For such a diagram, the Thuston--Bennequin number is the writhe minus the number of right cusps (maximum points with respect to projection onto the $x$--axis.  The Thurston--Bennequin number is the maximum value of this, taken over all diagrams satisfying the crossing criteria. 
 
  Lenhard Ng~\cite{ng1} has shown that basic bounds for this invariant are strict in the case of 2--bridge knots.  For all knots with 9 or fewer crossings, Ng~\cite{ng2}  has determined the value for both possible orientations.  In the same paper he has determined the Thurston Bennequin number for prime knots with 10 crossing,  except for $10_{132}$ for which the number is unknown with one of the two orientations.  (A knot and its mirror image might have different Thurston--Bennequin numbers.) 
    
  {\bf Update January 1, 2007} Ng~\cite{ng3} has determined the value for all but five 11 crossing knots.   {\bf  Update July 7,  2007} Ng,  in an update of~\cite{ng2},  has computed the values of the 6 knots of 11 or fewer crossing for which the Thurston-Bennequin number was previously uknown. 
  {\bf Update September,  2008} Litherland has computed the values for most unknown cases of 12 crossing knots. The remaining cases are: $12n_{41, 119, 120, 121, 145, 153,  199, 200, 243,  260, 282, 310, 322, 351, 362, 368, 377, 403, 414, 425, 475, 523, 549}$.
 \vskip.05in
  
  {\bf Update  June 15,  2012} Dynnikov and Prasolov~\cite{dp} have shown that TB($K$) + TB($-K)$ = $-$ArcIndex($K$).  From this   the computation of the Thurston-Bennequin numbers for 12 crossing knots follows.
    \subsection{Tunnel Number}  For every knot $K$ there is a collection of disjoint embedded arcs in $S^3$ with  endpoints on the knot such that the complement of the interior of a regular neighborhood of the union of the arcs and knot is a solid handlebody.  The minimum number of arcs required is the tunnel number of the knot.  As an example,  all two bridge knots have tunnel number one. The computation of tunnel  numbers of knots through 10 crossings was completed in~\cite{msy}.

 \subsection{Turaev Genus} The  Turaev genus  of a knot was first defined in \cite{DFKLS}. Here is a simple definition. Every knot diagram has smoothings of type A and B. To construct the A smoothing,  locally orient the arcs at each crossing so that the crossing is right handed. Smooth each crossing so that orientation is preserved. Similarly, to construct the B smoothing, smooth so that orientations are inconsistent. The two smoothings produce collections of circles in the plane, say SA and SB, with sA and sB circles, respectively. These two collections of circles are naturally cobordant via a cobordism of genus g = (2 + c - sa - sb)/2. The minimum of this genus over all diagrams for the knot is called the Turaev genus.
 We have the following result.
\vskip.1in
\noindent{\bf Theorem} {\it K is alternating if and only if the Turaev genus of K is 0.}
\vskip.1in
Lowrance~\cite{low} has proved   that the Turaev genus is an upper bound for the width of the Heegaard Floer knot homology, minus 1. A similar bound for the Khovanov width was found by Manturov in Minimal diagrams of classical and virtual links. See also~\cite{CKS}. In~\cite{DFKLS} it is proved that the Turaev genus is bounded above by the crossing number minus the span of Jones polynomial.   {\bf Update September 17, 2008:} Abe and Kishimoto~ \cite{ak} have shown the Turaev genus for all nonalternating knots under 12 crossings is 1, except for $11n_{95}$ and $11n_{118}$. For these two remaining knots, it might be either be 1 or 2. 
 \subsection{Unknotting Number}  The minimum number of changes of crossing required to convert a diagram of a given knot into a diagram for the unknot, taken over all possible diagrams of the knot, is called the unknotting number of the knot.  {\bf Update July 9, 2005:}  Brendan Owens~\cite{ow} has shown that the unknotting numbers of $9_{10}, 9_{13}, 9_{35},$  $9_{38},$ $10_{53},$  $10_{101}$, and $10_{120}$ are all 3, thus completing the tabulation for 9--crossing knots.   {\bf Update January 15, 2006:} For 10 crossing knots, all those with  unknotting number 1 have been identified by Ozsv\'ath-Szab\'o and Gordon-Lueke~\cite{gl,os}.

 \subsection{Unknotting Number 1}  {\bf Update January 15, 2006:} All knots of unknotting number 1 with 10 or fewer crossings have been determined.  The final examples were resolved by Ozsv\'ath-Szab\'o~\cite{os} and Gordon-Lueke~\cite{gl}.  An initial survey of 11-crossing knots has been done by Slavik Jablan and Radmila Sazdanovic, results of which are included in KnotInfo~\cite{ki}.
 
 {\bf Update February 6, 2009:}  During the fall, 2008, Josh Greene reported that new results of his, based on a combination of tools coming from Heegaard-Floer theory and Donaldson's original restrictions on the intresection forms of smooth manifolds, are sufficient to rule out unknotting number 1 for the remaining 100 cases for 11 alternating crossing knots. While that work was in progress, Slaven Jabuka independently announced the resolution of several cases of unknotting number 1, using rational Witt class invariants. Since then, Eric Staron has, independently, also ruled out unknotting number 1 for most of the remaining cases of 11 alternating crossing knots; details of that work, which uses Heegaard-Floer homology, are available directly from Eric, at the University of Texas. 
 
 {\bf Update August 15, 2018:} Lisa  Piccirillo~\cite{pic} has proved that the Conway knot, $11n_{34}$ is not smoothly slice. This implies it is of four-genus 1 and that its concordance order is 2 or more.
 
 \newcommand{\etalchar}[1]{$^{#1}$}


\end{document}